\documentclass[journal]{IEEEtran}

\usepackage{graphics} 
\usepackage{epsfig} 
\usepackage{mathptmx} 
\usepackage{times} 
\usepackage{amsmath} 
\usepackage{amssymb}  
\usepackage{psfrag}
\usepackage{textcomp}
\usepackage{amsthm}
\usepackage{setspace}
\usepackage{hyperref}
\usepackage{wasysym}
\usepackage{latexsym}
\usepackage{xfrac}
\usepackage{relsize}
\usepackage{pmat}
\usepackage{url}
\usepackage{comment}
\usepackage{pdfsync}

\usepackage[cal=cm]{mathalfa} 

\DeclareMathOperator{\spn}{span}
\DeclareMathOperator*{\argmin}{arg\;min}
\DeclareMathOperator{\co}{co}
\newcommand{\del}{\nabla}

\newcommand{\1}{\mathbf{1}}
\newcommand{\0}{\mathbf{0}}

\newcommand{\sumiv}{\sum_{i\in\mathcal{V}}}

\newcommand{\PXi}{\mathbf{P}_{\Xi}}

\newcommand{\Pa}{\mathbf{P}_{\!\!\mathcal{A}}}

\newcommand{\Bgo}[1]{\bar{B}_{#1}(\Gamma_0)}

 \theoremstyle{remark} 
\newtheorem{assump}{A\!\!}[section]

\newtheorem{claim}{Claim}[section]
\newtheorem{example}{Example}[section]
\newtheorem{theorem}{Theorem}[section]
\newtheorem{defn}{Definition}[section]
\newtheorem{rem}{Remark}[section]

\newtheorem{constr}{Construction}[section]

\begin{document}
\title{A Characterization of Semiglobal, Practical, Asymptotic Stability for Gain-Parametrized Systems}

\author{Karla~Kvaternik,~\IEEEmembership{Member,~IEEE,}
\thanks{K. Kvaternik is with the Department
of Mechanical and Aerospace  Engineering, Princeton University. e-mail: (karlak@princeton.edu). This research was supported by the Army Research Office Grant W911NF-14-1-0431.}}

\maketitle

\begin{abstract}
We consider a general class of nonlinear, constrained, discrete-time systems whose dynamics are parametrized by a set of gains. We define the \emph{semiglobal, practical, asymptotic stability} (SPAS) of compact sets for this class of systems, and we provide a Lyapunov characterization of such sets. A set $\mathcal{A}$ that is SPAS with respect to a given system need not be an attractor for that system. Relative to existing characterizations of similar qualitative behaviors, our SPAS theorem does not require the existence of an asymptotically stable attractor associated to a nominal counterpart of the given dynamics.
\end{abstract}

\section{Introduction}
\IEEEPARstart{W}{e} propose an explicit definition of \emph{semiglobal, practical, asymptotic stability} (SPAS) of compact sets for a class of nonlinear, constrained, discrete-time systems parametrized by a set of tunable gains. Our main contribution is a theorem that characterizes SPAS directly in terms of certain properties of a Lyapunov function.

In the literature, SPAS is most often associated  with dynamical systems affected by non-vanishing perturbations, and is typically characterized as a consequence of a  robustness property of the stability of their ``nominal'' counterparts. The application of existing SPAS theorems usually entails identifying an attractor $\mathcal{A}$ for a nominal system, and establishing its asymptotic stability. Then, it may be concluded that for any arbitrarily large set $\hat{\mathcal{B}}$ (contained in the basin of attraction $\mathcal{B}$ of $\mathcal{A}$ for the nominal system), and for any arbitrarily small set $\check{\mathcal{A}}\subset \hat{\mathcal{B}}$ containing $\mathcal{A}$, \emph{there exists} a (sufficiently small) non-vanishing perturbation of the nominal dynamics and a set $\tilde{\mathcal{A}}\subset\check{\mathcal{A}}$, such that $\tilde{\mathcal{A}}$ is asymptotically stable for the perturbed system, with a basin of attraction $\hat{\mathcal{B}}$. Examples of existing SPAS theorems include Theorem 17 in \cite{TeelCSM}, Corollary 1 in \cite{TeelCDC00a} and Theorem 10 in \cite{KelletSCL04}.

In this literature, the word ``semiglobal'' refers to the fact that $\hat{\mathcal{B}}$, the basin of attraction of the perturbed system, can be  arbitrarily large within $\mathcal{B}$, while ``practical'' typically refers to the fact that $\check{\mathcal{A}}$ can be arbitrarily small (so long as it contains $\mathcal{A}$), provided the perturbation of the nominal dynamics is sufficiently small. This literature provides a variety of characterizations of SPAS-like qualitative behaviors, and it is not always clear that these are equivalent.

Instead of characterizing the qualitative behavior of SPAS as a consequence of the asymptotic stability of a nominal attractor, we provide an explicit definition of SPAS sets, and characterize the SPAS property directly in terms of a set of conditions on a Lyapunov-like function. We thereby avoid the need to identify a nominal system and establish the asymptotic stability of its attractor.

We phrase our definition of SPAS in terms of a set of gains parametrizing the system dynamics. Such gains may for example model fixed step-sizes in iterative optimization algorithms or gradient estimation schemes. We say that a set $\mathcal{A}$ is SPAS for a given system if for \emph{some} neighborhood $B_{\rho}(\mathcal{A})$ of $\mathcal{A}$, and for every arbitrarily large set $B_{\sigma}(\mathcal{A})\supset B_{\rho}(\mathcal{A})$, \emph{there exists a set of gains} such that for all systems with gains in this set,  trajectories initialized inside $B_{\sigma}(\mathcal{A})$ are asymptotically attracted to $B_{\rho}(\mathcal{A})$ and never deviate far from $\mathcal{A}$ (q.v. $\S$\ref{sec:defs}, for a more precise definition).

Our definition thus retains the general intent of existing SPAS characterizations, while our SPAS theorem provides an alternative means for establishing the qualitative behavior that SPAS implies. This alternative characterization of SPAS is potentially useful in situations where a nominal version of a given system and its associated attractor may be difficult to identify or analyse (q.v. Examples \ref{ex:CO} and \ref{ex:cookedUp} in $\S$\ref{sec:examples}). In fact, a set $\mathcal{A} $ that satisfies the conditions of our SPAS theorem need not constitute a set of fixed points for the given dynamics.

For gain-parametrized systems having an obvious nominal counterpart with an asymptotically stable attractor, existing SPAS characterizations apply, and lead to conclusions that are consistent with our own; in such cases a Lyapunov function and attractor for the nominal system will often satisfy the conditions of our SPAS theorem, and the actual system dynamic can be expressed as a perturbed version of the nominal dynamic. However, whereas existing characterizations typically posit the existence of a non-vanishing perturbation corresponding to the desired size of the aforementioned sets $\hat{\mathcal{B}}$ and $\check{\mathcal{A}}$, our formulation is tantamount to assuming that a disturbance is \emph{given}, and that its size is not subject to design. Instead, our SPAS theorem posits the existence of a set of gains for which SPAS holds. In this sense, the SPAS characterization provided here is less abstract, for the class of gain-parametrized systems.

In the proof of our SPAS Theorem \ref{thm:spas}, we draw on some of the analytic techniques found in $\S$5.14 of \cite{Agarwal}, where a related notion of ``practical stability'' is considered. In particular, we model the invariance arguments used in the proof of Theorem \ref{thm:spas} on the proof of Theorem 4.14.2 in \cite{Agarwal}. Our SPAS theorem can be viewed as a generalization of Theorem 4.14.2 and Corollary 5.14.3 in \cite{Agarwal} to parametrized systems whose parameter valuations control the size of the basin of attraction and the size of the neighborhood to which trajectories initiated in this basin ultimately converge.

This paper is organized as follows. In $\S$\ref{sec:defs} we describe the class of systems under consideration, state a definition of SPAS sets for such systems, and we state a set of conditions characterizing the SPAS property in Theorem \ref{thm:spas}. In $\S$\ref{sec:examples} we demonstrate the utility and application of Theorem \ref{thm:spas} through a number of examples, and in $\S$\ref{sec:proof} we provide the proof of Theorem \ref{thm:spas}. We conclude the paper in $\S$\ref{sec:conclusions}.

\paragraph{Notation and Preliminaries}
Throughout this paper we use the Euclidean vector norm $\|\cdot\|$. If $S\subset\mathbb{R}^n$ is closed and $x_o\in\mathbb{R}^n$, then
\begin{equation}\label{eq:proj}
\mathbf{P}_S(x_o) = \argmin_{x\in S} \|x_o-x\|
\end{equation}
is the orthogonal projection of $x_o$ onto $S$. The set of non-negative real numbers is denoted by $\mathbb{R}_+$, while the set of positive real numbers is denoted by $\mathbb{R}_{++}$. For a point $x_o\in\mathbb{R}^n$, and $r\in\mathbb{R}_{++}$, $\bar{B}_r(x_o)=\{x\in\mathbb{R}^n: \|x-x_o\|\leq r\}$ and $B_r(x_o)=\{x\in\mathbb{R}^n: \|x-x_o\|<r\}$. For a compact set $S\subset\mathbb{R}^n$, $\bar{B}_{r}(S)=\{x\in\mathbb{R}^n\;\big|\; \|x-\mathbf{P}_S(x)\|\leq r\}$, while $B_{r}(S)=\{x\in\mathbb{R}^n\;\big|\; \|x-\mathbf{P}_S(x)\|< r\}$.
The set of continuous (resp. continuously differentiable) functions from $\mathbb{R}^n$ into $\mathbb{R}^m$ is denoted by $C^0[\mathbb{R}^n,\mathbb{R}^m]$ (resp. $C^1[\mathbb{R}^n,\mathbb{R}^m]$). We say that a function $V\in C^0[\mathbb{R}^n,\mathbb{R}_+]$ is \emph{positive definite} with respect to a closed set $S$ on a set $\Omega\supset S$ if, $V(S)=\{0\}$, and $V(x)>0$ for all $x\in\Omega\backslash S$. A function $V:x\mapsto V(x)$ on $\mathbb{R}^n$ is \emph{radially unbounded} with respect to a closed set $S\subset\mathbb{R}^n$, if for any $B\in\mathbb{R}$, there exists an $r\in\mathbb{R}_{++}$ such that $V(x)>B$, for all $x\in\mathbb{R}^n\backslash \bar{B}_r(S)$. The gradient of a differentiable function $J:\mathbb{R}^n\rightarrow\mathbb{R}$ is denoted by $\del J(\cdot)$.
For a sequence $(x(t))_{t=0}^{\infty}$ we use $x$ to stand for $x(t)$, and $x^+$ to stand for $x(t+1)$. For $S\subset\mathbb{R}^n$, $\co(S)$ is its convex hull.

\section{Semiglobal, Practical, Asymptotic Stability}\label{sec:defs}

We consider a class of discrete-time dynamical systems of the form
\begin{equation}\label{eq:system}
\xi^+ = \PXi\big[f(\xi;\pi)\big],\quad \xi\in\mathbb{R}^n,
\end{equation}
where $\Xi\subset\mathbb{R}^n$ is a closed, convex set, and $\pi\in\mathbb{R}^p$ parametrizes the function ${f:\mathbb{R}^n\rightarrow\mathbb{R}^n}$. 

We define \emph{semiglobal, practical, asymptotic stability} (SPAS) of compact sets for \eqref{eq:system} in Definitions \ref{def:stabilityNew2} to \ref{def:spasNew2}, and in Theorem \ref{thm:spas} we provide a Lyapunov characterization of SPAS. We discuss these definitions and characterization in the remarks that follow.

\begin{defn}\label{def:stabilityNew2}
A set $\mathcal{A}\subset\Xi$ is \emph{practically stable} for \eqref{eq:system} if for some $\check{\rho}_s\in\mathbb{R}_{++}$, and for any $\rho_s>\check{\rho}_s$, there exists a positive, real number $\delta$ and a set $P_s\subset\mathbb{R}^p$,  such that whenever $\pi\in P_s$ and $\xi(0)\in \bar{B}_{\delta}(\mathcal{A})\cap\Xi$, $\xi(t)\in\bar{B}_{\rho_s}(\mathcal{A})\cap\Xi$, for all $t\in\mathbb{N}$.
\end{defn}

\begin{defn}\label{def:attractivityGeneric}
A compact set $S\subset\mathbb{R}^n$ is  \emph{uniformly attractive for \eqref{eq:system} on a compact} $\Omega\subset\mathbb{R}^n$, if for every $\varepsilon\in\mathbb{R}_{++}$ for which $\bar{B}_{\varepsilon}(S)\cap\Xi\subset\Omega\cap\Xi$, there exists a number $T\in\mathbb{N}$ such that $\xi(t)\in\bar{B}_{\varepsilon}(S)$, whenever $\xi(0)\in\Omega$ and $t\geq T$.
\end{defn}

\begin{rem}
If a compact set $S\subset\mathbb{R}^n$ is  \emph{uniformly attractive} for \eqref{eq:system} on every compact $\Omega\subset\mathbb{R}^n$, then $S$ satisfies the usual definition of attractivity for \eqref{eq:system}.
\end{rem}

\begin{defn}\label{def:attractivityNew2}
A compact set $\mathcal{A}\subset\Xi$ is  \emph{semiglobally, practically attractive} for \eqref{eq:system} if for some $\check{\rho}_a\in\mathbb{R}_{++}$, and for any $\sigma,\rho_a,\in\mathbb{R}_{++}$, with $\sigma>\rho_a>\check{\rho}_a$, there exists a set $P_a\subset\mathbb{R}^p$ such that whenever $\pi\in P_a$, the set $\bar{B}_{\rho_a}(\mathcal{A})$ is uniformly attractive for \eqref{eq:system} on $\bar{B}_{\sigma}(\mathcal{A})$.
\end{defn}


\begin{defn}\label{def:spasNew2}
A set $\mathcal{A}\subset\Xi$ is \emph{semiglobally practically asymptotically stable} (SPAS)  for \eqref{eq:system} if it is practically stable and semiglobally, practically attractive for \eqref{eq:system}.
\end{defn}

The following theorem characterizes SPAS in terms of a Lyapunov function with certain properties.

\begin{theorem}\label{thm:spas}
Consider the system \eqref{eq:system}, and suppose there exists a function $V\in C^0[\mathbb{R}^n, \mathbb{R}_+]$ which is radially unbounded and positive definite with respect to a compact set $\mathcal{A}\subset\Xi$ on $\mathbb{R}^n$. Suppose that for some $\epsilon_o\in\mathbb{R}_+$ and for any positive, real $\sigma_o$, $\rho_o$ and $b_o$ (with $\sigma_o>\epsilon_o+\rho_o$) 
there exists a set $P_o\subset\mathbb{R}^p$ and a function $W_{\sigma_o,\epsilon_o}\in C^0[\mathbb{R}^n,\mathbb{R}]$ such that whenever $\pi\in P_o$:
\begin{itemize}
\item P1: $W_{\sigma_o,\epsilon_o}(\xi) > 0$ for all  $\xi\in \Xi\cap (\bar{B}_{\sigma_o}(\mathcal{A})\backslash B_{\epsilon_o+\rho_o}(\mathcal{A}))$,
\item P2: $\Delta V(\xi)\leq -W_{\sigma_o,\epsilon_o}(\xi)$, for all $\xi\in\Xi\cap (\bar{B}_{\sigma_o}(\mathcal{A})\backslash B_{\epsilon_o+\rho_o}(\mathcal{A}))$, and
\item P3: $\Delta V(\xi) \leq b_o$, for all $\xi\in \Xi\cap\bar{B}_{\epsilon_o+\rho_o}(\mathcal{A})$.
\end{itemize}

Then, $\mathcal{A}$ is SPAS for \eqref{eq:system}, with Lyapunov function $V(\cdot)$.
$\square$\end{theorem}

The proof is given in $\S$\ref{sec:proof}.

\begin{rem}
To satisfy the SPAS definition, the numbers $\check{\rho}_s$ and $\check{\rho}_a$ in Definitions \ref{def:stabilityNew2} and \ref{def:attractivityNew2} need not coincide, nor do the sets $P_s$ and $P_a$.
$\diamondsuit$\end{rem}

\begin{rem}\label{rem:AnotAttractor}
As shown in the upcoming Examples \ref{ex:cookedUp} and \ref{ex:CO}, the set $\mathcal{A}$ need not be an attractor for \eqref{eq:system}.
$\diamondsuit$\end{rem}

\begin{rem}[$\epsilon_o$ determines $\check{\rho}_s$ and $\check{\rho}_a$]
The numbers $\check{\rho}_s$ and $\check{\rho}_a$ in Definitions \ref{def:stabilityNew2} and \ref{def:attractivityNew2} are determined by the number $\epsilon_o$ in the statement of Theorem \ref{thm:spas}. More specifically, in the proof of Theorem \ref{thm:spas}, we construct $\check{\rho}_a$ as $\check{\rho}_a=\epsilon_o$, and $\check{\rho}_s$ as increasing with $\epsilon_o$ and depending additionally only on the properties of the function $V(\cdot)$.

Depending on the application, $\epsilon_o$ may relate to the size of a \emph{given} non-vanishing disturbance whose effect cannot be affected by any of the gain parameters. Examples of systems for which $\epsilon_o>0$ are considered in \cite{KvaternikSPSP}.

In contrast to $\epsilon_o$, the number $\rho_o$ in the statement of  Theorem \ref{thm:spas} relates to those components of errors, typically arising within the Lyapunov analysis itself, whose size can be modulated through some of the system gains (q.v. Example \ref{ex:cookedUp}).
$\diamondsuit$\end{rem}


\begin{rem}[Relating a conventional definition of stability to Definition \ref{def:stabilityNew2}]
A conventional definition of set stability would apply to a particular $\pi_o$-instance of \eqref{eq:system}, and might read as follows:
\begin{defn}\label{def:conventionalStability}
A compact set $\mathcal{A}\subset\Xi$ is  \emph{stable} for a system $\xi^+ = \PXi\big[f(\xi;\pi_o)\big]$ if for every $\varepsilon\in\mathbb{R}_{++}$ there exists a $\delta\in\mathbb{R}_{++}$ such that $\xi(t)\in\bar{B}_{\varepsilon}(\mathcal{A})$ for all $t\in\mathbb{N}$, whenever $\xi(0)\in\bar{B}_{\delta}(\mathcal{A})$.
\end{defn}
Definition \ref{def:stabilityNew2}, which pertains to a parametrized family of dynamical systems, cannot be directly compared to Definition \ref{def:conventionalStability}, even for the case in which $\check{\rho}_s=0$. Specifically, for any given $\rho_s$ (analogous to $\varepsilon$ in Definition \ref{def:conventionalStability}), Definition \ref{def:stabilityNew2} requires the existence of a $\delta$, and only a sub-family of \eqref{eq:system} to generate trajectories that never exit $\bar{B}_{\rho_s}(\mathcal{A})$ when initialized within $\bar{B}_{\delta}(\mathcal{A})$. In other words, Definition \ref{def:stabilityNew2} allows the instantiations of \eqref{eq:system} to depend on the given $\rho_s$.

However, in terms of the conditions of Theorem \ref{thm:spas} (and based on the analytic techniques employed in the invariance arguments in its proof) conventional stability is obtained for the special case in which P1 to P3 hold with  $\epsilon_o=\rho_o= 0$, and independently of $\pi$. An example of a system for which this happens is given in Example 3.3 in \cite{KvaternikSPSP}.
%
$\diamondsuit$\end{rem}

\begin{rem}[Related notions of SPAS]
There is a variety of characterizations of SPAS in the literature. For example, Corollary 1 of \cite{TeelCDC00a}, which applies to continuous-time differential inclusions, can be paraphrased as follows: if the compact set $\mathcal{A}$ is asymptotically stable with basin of attraction $\mathcal{G}$ for the nominal system $\dot{x}\in F(x)$, and the multifunction $F(\cdot)$ satisfies certain technical conditions, then, for any arbitrarily large set $B\subset\mathcal{G}$ and any arbitrarily small set $A$ (contained in $\mathcal{G}$, and containing $\mathcal{A}$), \emph{there exists} a non-vanishing perturbation of size $\epsilon\in\mathbb{R}_{++}$ and a compact set
\begin{equation}\label{eq:A}
\mathcal{A}_{\epsilon}\subset A
\end{equation}
such that $\mathcal{A}_{\epsilon}$ is asymptotically stable for the perturbed system
\begin{equation*}
\dot{x}\in  B_{\epsilon}\big( \co F(B_{\epsilon}(x))\big),
\end{equation*}
with a basin of attraction containing $B$.

Another example is  Theorem 17 in \cite{TeelCSM}, which applies to a very general class of hybrid systems. This theorem can be paraphrased (and specialized to difference inclusions only) as follows: if the compact set $\mathcal{A}$ is asymptotically stable with basin of attraction $\mathcal{G}$ for the nominal system $x^+\in \mathbf{P}_{\Xi}[F(x)]$, and the multifunction $F(\cdot)$ satisfies certain technical conditions, then, for any arbitrarily small $\rho\in\mathbb{R}_{++}$ and any arbitrarily large set $B\subset\mathcal{G}$, there exists a non-vanishing perturbation of size $\epsilon\in\mathbb{R}_{++}$, a proper indicator function  $\omega:\mathbb{R}^n\rightarrow\mathbb{R}_{++}$ for $\mathcal{A}$ on $\mathbb{R}^n$ (q.v. p53 in \cite{TeelCSM} for a definition), and  a $\mathcal{KL}$-class function  $\beta(\cdot,\cdot)$ such that all solutions of the perturbed system
\begin{equation}
x^+\in B_{\epsilon}\big(\mathbf{P}_{\Xi}[F(B_{\epsilon}(x))]),
\end{equation}
initialized inside $B$, satisfy
\begin{equation}\label{eq:klfunct}
\omega(x(t))\leq \beta(\omega(x(0)),t) + \rho, \quad \forall t\in\mathbb{N}.
\end{equation}

One difference between characterizations such as these and the characterization provided in Definition \ref{def:spasNew2} and Theorem \ref{thm:spas}, is that the latter does not require the existence of an asymptotically stable attractor for a nominal counterpart to \eqref{eq:system}. This alternative characterization could therefore be useful in the analysis of systems for which a nominal counterpart and its associated attractor may be difficult to identify and analyse, and for which it is easier to directly identify a set $\mathcal{A}$ that might approximate the attractor of the given system. As noted in Remark \ref{rem:AnotAttractor}, the set $\mathcal{A}$ in the statement of Theorem \ref{thm:spas} need not be an actual attractor for \eqref{eq:system}. In Example \ref{ex:cookedUp} we consider an iterative numerical method that exemplifies this situation.

For the class of gain-parametrized systems \eqref{eq:system} for which a nominal dynamic and an asymptotically stable nominal attractor can be easily identified,  existing SPAS characterizations do apply, but lead to conclusions that are different, though consistent with those of Theorem \ref{thm:spas}. Specifically, whereas existing characterizations of SPAS typically posit the existence of a non-vanishing perturbation  corresponding to the desired size of the sets $B$ and $A$ in \eqref{eq:A}, or the number $\rho$ in \eqref{eq:klfunct},  we assume that a disturbance is \emph{given} and that its size, which relates to $\epsilon_o$ in Theorem \ref{thm:spas}, is not subject to design. Instead, Theorem \ref{thm:spas} posits the existence of a set of gains for which \eqref{eq:system} exhibits the qualitative behavior specified by Definition \ref{def:spasNew2}. In this sense, for the class of gain-parametrized systems represented by \eqref{eq:system}, the SPAS characterization that we provide is less abstract than existing SPAS characterizations.
$\diamondsuit$\end{rem}


\section{Examples}\label{sec:examples}
In the following three examples, we motivate Theorem \ref{thm:spas} and demonstrate its application.

\begin{example}[Generic iterative numerical methods]\label{ex:cookedUp}
Consider an iterative numerical method of the form
\begin{equation}\label{eq:basicSys}
y^+ = \mathbf{P}_{\Xi}(y-\alpha s(y)),
\end{equation}
where $\Xi$ is a closed, convex subset of $\mathbb{R}^n$, $\alpha\in\mathbb{R}_{++}$ is a tunable gain, $y\in\mathbb{R}^n$, and $s(\cdot)$ is locally Lipschitz continuous.

Suppose we know that for some compact set $\mathcal{A}\subset\mathbb{R}^n$ and number $\tau\in\mathbb{R}_{++}$, the following hold:
\begin{align}
& (y-\Pa (y))^Ts(y)\geq \tau \|y-\Pa (y)\|^2,\quad \forall y\in\mathbb{R}^n\backslash\mathcal{A},\label{eq:ineq}\\
& \exists y\in \mathcal{A}\;|\; s(y)\neq 0. \label{eq:nonEq}
\end{align}
This can happen, for example, when $\mathcal{A}$ is a set of minima of some strongly convex function, and $s(\cdot)$ is an approximation of its gradient. We can then use Theorem \ref{thm:spas} to conclude that $\mathcal{A}$ is SPAS for \eqref{eq:basicSys} despite the fact that $\mathcal{A}$ does not constitute a set of fixed points for \eqref{eq:basicSys}, and therefore cannot be an attractor for \eqref{eq:basicSys}.

Specifically, we let $V(y) = \tfrac{1}{2}\|y-\Pa (y)\|^2$ and $\Delta V(y) = V(y^+)-V(y)$. By expanding the expression for $V(\cdot)$ and using the properties of the projection operator (q.v. Lemma 4.1, \cite{KvaternikSPSP}), it can be shown that
\begin{equation*}
\Delta V(y) \leq -\alpha (y-\Pa (y))^Ts(y) + \tfrac{1}{2}\alpha^2 \|s(y)\|^2.
\end{equation*}
Using \eqref{eq:ineq}, Young's inequality and the fact that $s(\cdot)$ is locally Lipschitz continuous, we have that for any $\sigma_o\in\mathbb{R}_{++}$, there exists a number $L_s\in\mathbb{R}_{++}$ such that
\begin{align}
\Delta V(y) &\leq -\alpha\tau\|y-\Pa (y)\|^2 + \alpha^2L_s^2 \|y-\Pa(y)\|^2
\nonumber \\
&\qquad \quad+ \alpha^2\|s(\Pa(y))\|^2,
\end{align}
for all $y\in\bar{B}_{\sigma_o}(\mathcal{A})$.

Since $s(\mathcal{A})\neq \{0\}$ and $\mathcal{A}$ is compact, the number
\begin{equation}
s^* = \max_{y\in\mathcal{A}} \|s(\Pa (y))\|^2
\end{equation}
exists, and is positive. Therefore, for all $y\in\bar{B}_{\sigma_o}(\mathcal{A})$,
\begin{align}
\Delta V(y) &\leq -\alpha (\tau - \alpha L_s^2) \|y-\Pa (y)\|^2 + \alpha^2s^*,
\end{align}
and we see that for any given $\rho_o$, $\sigma_o$ and $b_o$, the conditions of Theorem \ref{thm:spas} are satisfied with $\epsilon_o=0$ and
\begin{equation}
W_{\sigma_o,0} = \alpha (\tau - \alpha L_s^2) \big( \|y-\Pa (y)\|^2 - \tfrac{\alpha s^*}{\tau - \alpha L_s^2}\big),
\end{equation}
by taking $P_o=(0,\hat{\alpha}]$, where $\hat{\alpha}=\min\{\alpha_{b_o},\alpha_{\rho_o},\alpha_W\}$, and
\begin{align*}
\alpha_{b_o} &= \sqrt{\tfrac{b_o}{s^*}}
\\
\alpha_{\rho_o} &= \tfrac{\rho_o^2\tau}{s^*+\rho_o^2L_s^2}
\\
\alpha_W & = \tfrac{\tau}{L_s^2}.
\end{align*}
$\diamondsuit$\end{example}

\begin{example}[Iterative methods with SPSP search directions]
Another motivation for the way in which Theorem \ref{thm:spas} is formulated is that it is particularly well suited to the analysis of a large class of iterative methods employing \emph{semiglobal, practical, strictly pseudogradient} (SPSP) search directions. This class of systems takes the form
\begin{equation}\label{eq:basicSys2}
y^+ = \mathbf{P}_{\Xi}(y-\alpha s),
\end{equation}
where at every $y\in\mathbb{R}^n$, the \emph{search direction} $s$ can be any element of a multifunction $\Psi:\mathbb{R}^n\rightrightarrows \mathbb{R}^n$ that has the SPSP property, which is defined as follows \cite{KvaternikSPSP}:

\emph{A multifunction
$\Psi:\mathbb{R}^n\rightrightarrows \mathbb{R}^n$ is SPSP on a set $\Xi\subset\mathbb{R}^n$ with respect to a differentiable function $V:\mathbb{R}^n\rightarrow\mathbb{R}_+$, which is positive definite and radially unbounded with respect to a compact set $\mathcal{A}\subset\mathbb{R}^n$, if for some $\epsilon\in\mathbb{R}_+$ and $b\in\mathbb{R}_+$, and for any $\sigma\in\mathbb{R}_{++}$  (with $\sigma >\epsilon$),
\begin{equation*}
\del V(y)^Ts\geq -b,\quad \forall y\in\Xi\cap\bar{B}_{\epsilon}(\mathcal{A}),\;\forall s\in\Psi(y),
\end{equation*}
and there exists a function $\phi_{\sigma,\epsilon}\in C^0[\mathbb{R}^n,\mathbb{R}]$ which is positive on $\Xi\cap \big(\bar{B}_{\sigma}(\mathcal{A})\backslash B_{\epsilon}(\mathcal{A})\big)$ and radially unbounded with respect to $\mathcal{A}$ on $\Xi$, such that
\begin{equation*}
\del V(y)^Ts\geq \phi_{\sigma,\epsilon}(y), \quad \forall y\in\Xi\cap (\bar{B}_{\sigma}(\mathcal{A})\backslash B_{\epsilon}(\mathcal{A})),\;\forall s\in\Psi(y).
\end{equation*}
}

In \cite{KvaternikSPSP}, we consider a number of optimization-related examples of iterative methods that fall within this class, we show that under a variety of standard additional assumptions on the search directions, the set $\mathcal{A}$ is SPAS for this class of algorithms, and we prove that the SPSP property is robust under absolute and relative deterministic errors on the search directions.
$\diamondsuit$\end{example}

\begin{example}[Consensus Optimization]\label{ex:CO}
Consider a system of the form
\begin{align}\label{eq:co}
x_i^+ &= \sumiv [A]_{i,j} x_j - \alpha \del J_i(x_i),\quad i\in\mathcal{V},
\end{align}
where $\mathcal{V}=\{1,\ldots,N\}$, $\forall i\in\mathcal{V}$ $x_i\in\mathbb{R}$ and  $J_i:\mathbb{R}\rightarrow\mathbb{R}$ is differentiable and strictly convex, and $A\in\mathbb{R}^{N\times N}$ is stochastic, symmetric, and such that there exists a number $\mu\in [0,1)$ for which
\begin{equation}\label{eq:mu}
\|A z\|^2 \leq \mu \|z\|^2,
\end{equation}
whenever $z\in\spn\{\1_N\}^{\bot}$, where $\1_N=[1,\ldots,1]^T\in\mathbb{R}^N$. The strict convexity of $J_i(\cdot)$ implies the existence of a unique minimizer $x_i^*\in\mathbb{R}$. Algorithm \eqref{eq:co} is a special case of a decentralized optimization method originating in \cite{NedicTAC09}, and further studied in \cite{KvaternikPHD}.

Unless it happens to be the case that the individual optima coincide,  the actual equilibrium of \eqref{eq:co} may be difficult to identify. An equilibrium of \eqref{eq:co} is a fixed point of the mapping
\begin{equation}\label{eq:mapping}
F(\mathbf{x}) := A\mathbf{x}-\alpha \mathbf{s}(\mathbf{x}),
\end{equation}
where $\mathbf{x}=[x_1,\ldots, x_N]^T$, and $\mathbf{s}(\mathbf{x}) = [\del J_1(x_1),\ldots,\del J_N(x_N)]^T$. From this expression it is evident that a point  $\mathbf{x}_o:=[x_1^*,\ldots,x_N^*]^T$, satisfying $\mathbf{s}(\mathbf{x}_o)=\0_N$, is not an equilibrium for \eqref{eq:co} unless $\mathbf{x}_o\in\spn\{\1_N\}$; the properties of $A$ imply that  $A\mathbf{x} = \mathbf{x}$ only for $\mathbf{x}\in\spn\{\1_N\}$. On the other hand for $\mathbf{x}\in\spn\{\1_N\}$, it is not generally the case that $\mathbf{s}(\mathbf{x})=\0_N$.

One way to analyse this system is to resolve its dynamics along the so-called ``agreement subspace'' $\spn\{\1_N\}$ and its orthogonal complement. To that end, it can be shown (see $\S$3.4.1 in \cite{KvaternikPHD}) that the variables $y=\tfrac{1}{N}\1_M^T x$ and $z=M x$, where  $M=\1_N-\tfrac{1}{N}\1_N\1_N^T$, evolve according to
\begin{align}\label{eq:yev}
y^+ &= y - \tfrac{\alpha}{N}\1_N^T \mathbf{s}(z+\1_Ny)
\end{align}
and
\begin{equation}\label{eq:zev}
z^+ =  A z - \alpha M \mathbf{s}(z+\1_Ny).
\end{equation}
We observe that $\1_N^T\mathbf{s}(\1_Ny)=\del J(y)$, where $J(\cdot)$ denotes the sum of individual objectives $J_i(\cdot)$, and we note that if $J_i(\cdot)$ is strictly convex for each $i\in\mathcal{V}$, then $J(\cdot)$ is strictly convex, and has a unique minimizer $x^*\in\mathbb{R}$.

Therefore, when $z\equiv 0$, \eqref{eq:yev} takes the form of \eqref{eq:basicSys}, with $\Xi=\mathbb{R}^N$ and $s(y)=\tfrac{1}{N} \del J(y)$. Moreover, $s(\cdot)$ has the \emph{strict pseudogradient}  property (a special case of the SPSP property, defined in \cite{KvaternikSPSP}) with respect to either $\|y-x^*\|^2$ or $J(y)$, as shown in Example 3.1 in \cite{KvaternikSPSP}. This observation, together with \eqref{eq:mu}, is exploited in $\S$3.6 of \cite{KvaternikPHD} to show that the function $V(y,z)=N\|y-x^*\|^2+\|z\|^2$  satisfies the conditions of Theorem \ref{thm:spas}. We are thereby able to conclude that the set $\mathcal{A}=\{x^*\}\times\{\0_N\}$, though not generally an equilibrium, is SPAS for  \eqref{eq:yev}-\eqref{eq:zev}.
$\diamondsuit$\end{example}

\section{Proof of Theorem \ref{thm:spas}}\label{sec:proof}

In order to satisfy Definition \ref{def:spasNew2}, we must demonstrate that given the hypotheses of Theorem \ref{thm:spas} and any numbers $\sigma$, $\rho_s$ and $\rho_a$, there exist parameter sets $P_s$ and $P_a$, and numbers $\check{\rho}_s$ and $\check{\rho}_a$ for which the behavior of \eqref{eq:system} specified in Definitions \ref{def:stabilityNew2} and \ref{def:attractivityNew2} is guaranteed.

\subsubsection{Preliminaries}

For any $r\in\mathbb{R}_{+}$, we let $\Gamma_r$ denote the set $\{\xi\in\mathbb{R}^n\;\big|\; V(\xi)\leq r\}$ while $\partial \Gamma_r$ denotes its boundary $\{\xi\in\mathbb{R}^n\;\big|\; V(\xi)= r\}$.
We observe that since $V(\cdot)$ is assumed to be radially unbounded with respect to $\mathcal{A}$, and $\mathcal{A}\equiv \Gamma_0$ is assumed to be compact, all the sublevel sets $\Gamma_r$ are compact. Moreover, for any $r\in\mathbb{R}_{++}$, the intersection $\Gamma_r\cap\Xi$ is nonempty since $\mathcal{A}$ is a subset of both $\Gamma_r$ and $\Xi$.

We begin with a set construction and a claim which are applied in the proofs of both practical stability and semiglobal, practical attractivity.

\begin{figure}[htb!]
      \centering
      \includegraphics[scale=.3]{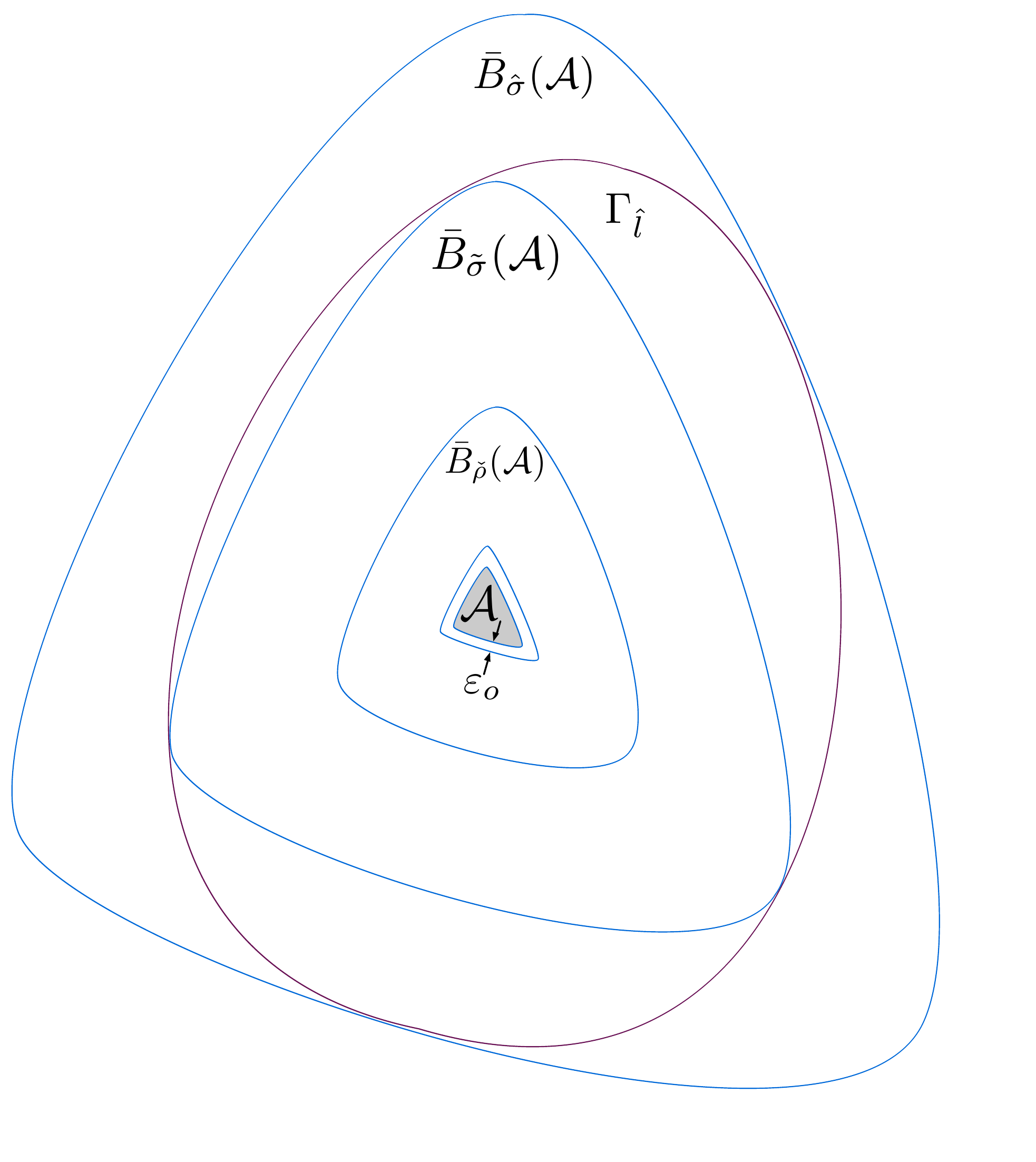}
      \caption{An illustration of Construction \ref{con:outer}, used in the proof of Theorem \ref{thm:spas}.}
      \label{fig:outer}
      \vspace*{-.1in}
   \end{figure}

\begin{constr}[Picking $\hat{\sigma}$ and $\hat{l}$ from $\tilde{\sigma}$]\label{con:outer}
Let $\tilde{\sigma}>\epsilon_o$ be an arbitrarily large, positive, real number, and let $\Gamma_{\hat{l}}$ be the smallest sublevel set of $V(\cdot)$ containing $\bar{B}_{\tilde{\sigma}}(\mathcal{A})$ -- i.e.,
\begin{align}\label{eq:lHat}
\hat{l} &= \max_{\xi\in\mathbb{R}^n}  V(\xi), \quad \text{s.t.} \quad  \|\xi-\mathbf{P}_{\!\mathcal{A}}(\xi)\| = \tilde{\sigma}.
\end{align}

Let $\bar{B}_{\hat{\sigma}}(\mathcal{A})$ be the smallest ball containing $\Gamma_{\hat{l}}$ -- i.e.,
\begin{align}\label{eq:sigma_o}
\hat{\sigma} &= \max_{\xi\in\partial \Gamma_{\hat{l}}}  \|\xi - \mathbf{P}_{\!\mathcal{A}}(\xi)\|.
\end{align}
$\diamondsuit$\end{constr}
Construction \ref{con:outer} is illustrated in Figure \ref{fig:outer}

The following claim states that certain sublevel sets of $V(\cdot)$ can be made forward invariant for \eqref{eq:system} under an appropriate restriction on the  gains $\pi$.

\begin{claim}[Forward invariance of $\Gamma_l$]\label{cl:invarianceNew}
Let $\tilde{\sigma}$ and $\tilde{\rho}$ be any positive, real numbers such that $\bar{B}_{\epsilon_o+\tilde{\rho}}(\mathcal{A})\subset \bar{B}_{\tilde{\sigma}}(\mathcal{A})$, and use Construction \ref{con:outer} to generate the numbers $\hat{\sigma}$ and $\hat{l}$ from $\tilde{\sigma}$. Then, for any $l\in (\tilde{l},\hat{l}]$, with $\tilde{l}$ constructed so that $\Gamma_{\tilde{l}}$ is the smallest sublevel set of $V(\cdot)$ containing $\bar{B}_{\epsilon_o+\tilde{\rho}}(\mathcal{A})$ -- i.e.,
\begin{align}\label{eq:lo}
\tilde{l} &= \max_{\xi\in\mathbb{R}^n}  V(\xi), \quad \text{s.t.} \quad  \|\xi-\mathbf{P}_{\!\mathcal{A}}(\xi)\| = \epsilon_o+\tilde{\rho},
\end{align}
there exists a parameter set $P_l$ such that whenever $\pi\in P_l$, $\Gamma_l$ is forward invariant for \eqref{eq:system}.
\proof
By the hypotheses of Theorem \ref{thm:spas}, there exists a parameter set $P_l$, corresponding to the choice $(\sigma_o,\rho_o,b_o)=(\hat{\sigma},\tilde{\rho},l-\tilde{l})$, such that
\begin{align}
&\Delta V(\xi) \leq l-\tilde{l} \quad \forall \xi\in  \bar{B}_{\epsilon_o+\tilde{\rho}}(\mathcal{A})\cap\Xi,\quad\text{and}
\label{eq:DeltaVbo}
\\
& \Delta V(\xi) < 0, \quad \forall \xi\in \big(\bar{B}_{\hat{\sigma}} \backslash B_{\epsilon_o+\tilde{\rho}}(\mathcal{A})\big)\cap\Xi,
\label{eq:DeltaVneg}
\end{align}
whenever $\pi\in P_l$.

Let $\pi\in P_l$. We will show that if for some $t\in\mathbb{N}$, $\xi\in \Gamma_l\cap\Xi$, then necessarily $\xi^+\in \Gamma_l\cap\Xi$. Suppose that $\xi\in \Gamma_l\cap\Xi$. Then, either $\xi\in \bar{B}_{\epsilon_o+\tilde{\rho}}(\mathcal{A})\cap\Xi$, or $\xi\in (\Gamma_l\backslash B_{\epsilon_o+\tilde{\rho}}(\mathcal{A}))\cap\Xi$.

If $\xi\in \bar{B}_{\epsilon_o+\tilde{\rho}}(\mathcal{A})\cap\Xi$, then $\xi\in \Gamma_{\tilde{l}}\cap\Xi$, meaning that $V(\xi)<\tilde{l}$. By \eqref{eq:DeltaVbo}, $V(\xi^+)-V(\xi)\leq l-\tilde{l}$, and therefore $V(\xi^+)\leq l$, meaning that $\xi^+\in \Gamma_l\cap\Xi$.

On the other hand if $\xi\in (\Gamma_l\backslash B_{\epsilon_o+\tilde{\rho}}(\mathcal{A}))\cap\Xi$, then $V(\xi)\leq l$, while $\Delta V(\xi)<0$, by \eqref{eq:DeltaVneg}. Therefore, $V(\xi^+)\leq l$, implying that $\xi^+\in \Gamma_l\cap\Xi$. The conclusion then follows from the principle of induction.
$\diamondsuit$\end{claim}

   \begin{figure}[htb!]
      \centering
      \includegraphics[scale=.6]{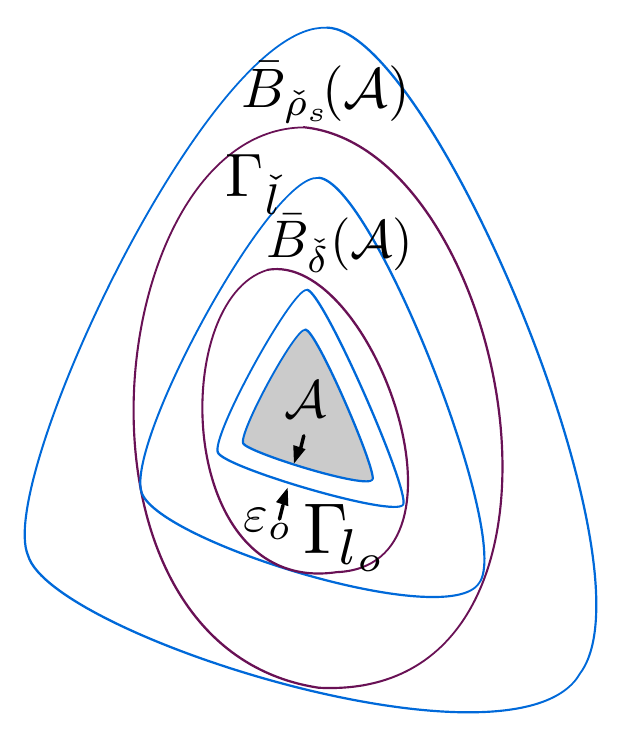}
      \caption{An illustration of Construction \ref{con:inner}, used in the proof of Theorem \ref{thm:spas}.
      }
      \label{fig:inner}
      \vspace*{-.1in}
   \end{figure}

\subsubsection{Practical Stability}\label{sec:stability}

To show that \eqref{eq:system} exhibits practical stability at $\mathcal{A}$, we first construct the requisite number $\check{\rho}_s$, and then apply Construction \ref{con:backward} and Claim \ref{cl:invarianceNew} to generate the parameter set $P_s$.

\begin{constr}[Constructing $\check{\rho}_s$ from $\epsilon_o$]\label{con:inner}
Let $\epsilon_o\in\mathbb{R}_{+}$ be as in the theorem statement, and let $\Gamma_{l_o}$ be the smallest sublevel set of $V(\cdot)$ containing $\bar{B}_{\epsilon_o}(\mathcal{A})$ -- i.e.,
\begin{align}\label{eq:lo}
l_o &= \max_{\xi\in\mathbb{R}^n}  V(\xi), \quad \text{s.t.} \quad  \|\xi-\mathbf{P}_{\!\mathcal{A}}(\xi)\| = \epsilon_o.
\end{align}

Let $\bar{B}_{\check{\delta}}(\mathcal{A})$ be the smallest ball containing $\Gamma_{l_o}$ -- i.e.,
\begin{align}\label{eq:deltaCheck}
\check{\delta} &= \max_{\xi\in\partial \Gamma_{l_o}}  \|\xi - \mathbf{P}_{\!\mathcal{A}}(\xi)\|.
\end{align}

Let $\Gamma_{\check{l}}$ be the smallest sublevel set of $V(\cdot)$ containing $\bar{B}_{\check{\delta}}(\mathcal{A})$ -- i.e.,
\begin{align}\label{eq:lCheck}
\check{l} &= \max_{\xi\in\mathbb{R}^n}  V(\xi), \quad \text{s.t.} \quad  \|\xi-\mathbf{P}_{\!\mathcal{A}}(\xi)\| = \check{\delta}.
\end{align}

Let $\bar{B}_{\check{\rho}_s}(\mathcal{A})$ be the smallest ball containing $\Gamma_{\check{l}}$ -- i.e.,
\begin{align}\label{eq:rhoCheck}
\check{\rho}_s &= \max_{\xi\in\partial \Gamma_{\check{l}}}  \|\xi - \mathbf{P}_{\!\mathcal{A}}(\xi)\|.
\end{align}
$\diamondsuit$\end{constr}
Construction \ref{con:inner} is depicted in Figure \ref{fig:inner}.

We use the next construction to generate the number $\delta$ in Definition \ref{def:stabilityNew2} from the given $\rho_s$. We also construct the numbers $\rho_{o,s}$ and $l_{\rho_s}$, which are used as inputs to the Claim \ref{cl:invarianceNew} to generate a parameter set $P_s$.

\begin{constr}[Picking $\delta$, $\rho_{o,s}$ and $l_{\rho_s}$ from $\rho_s$]\label{con:backward}
Use Construction \ref{con:inner} to generate the number $\check{\rho}_s$, and let $\rho_s\in (\check{\rho}_s,\sigma)$ be arbitrary, as in Definition \ref{def:stabilityNew2}.

Let $\Gamma_{l_{\rho_s}}$ be the largest sublevel set of $V(\cdot)$ that is contained inside $\bar{B}_{\rho_s}(\mathcal{A})$ -- i.e.,
\begin{align}\label{eq:lRho}
l_{\rho_s} &= \min_{\xi\in\mathbb{R}^n}  V(\xi), \quad \text{s.t.} \quad  \|\xi-\mathbf{P}_{\!\!\mathcal{A}}(\xi)\| = \rho_s.
\end{align}

Let $\bar{B}_{\delta}(\mathcal{A})$ be the largest ball contained inside $\Gamma_{l_{\rho_s}}$ -- i.e.,
\begin{align}\label{eq:delta}
\delta &= \min_{\xi\in\partial \Gamma_{l_{\rho_s}}}  \|\xi - \mathbf{P}_{\!\!\mathcal{A}}(\xi)\|.
\end{align}

Let $\Gamma_{l_{\delta}}$ be the largest sublevel set of $V(\cdot)$ that is contained inside $\bar{B}_{\delta}(\mathcal{A})$ -- i.e.,
\begin{align}\label{eq:ldelta}
l_{\delta} &= \min_{\xi\in\mathbb{R}^n}  V(\xi), \quad \text{s.t.} \quad  \|\xi-\mathbf{P}_{\!\!\mathcal{A}}(\xi)\| = \delta.
\end{align}

Let $\bar{B}_{\epsilon_o+\rho_{o,s}}(\mathcal{A})$ be the largest ball contained inside $\Gamma_{l_{\delta}}$ -- i.e.,
\begin{align}\label{eq:rho_o}
\rho_{o,s} &= \min_{\xi\in\partial \Gamma_{l_{\delta}}}  \|\xi - \mathbf{P}_{\!\mathcal{A}}(\xi)\| - \epsilon_o.
\end{align}
%
$\diamondsuit$\end{constr}

\begin{figure}[htb!]
      \centering
      \includegraphics[scale=.5]{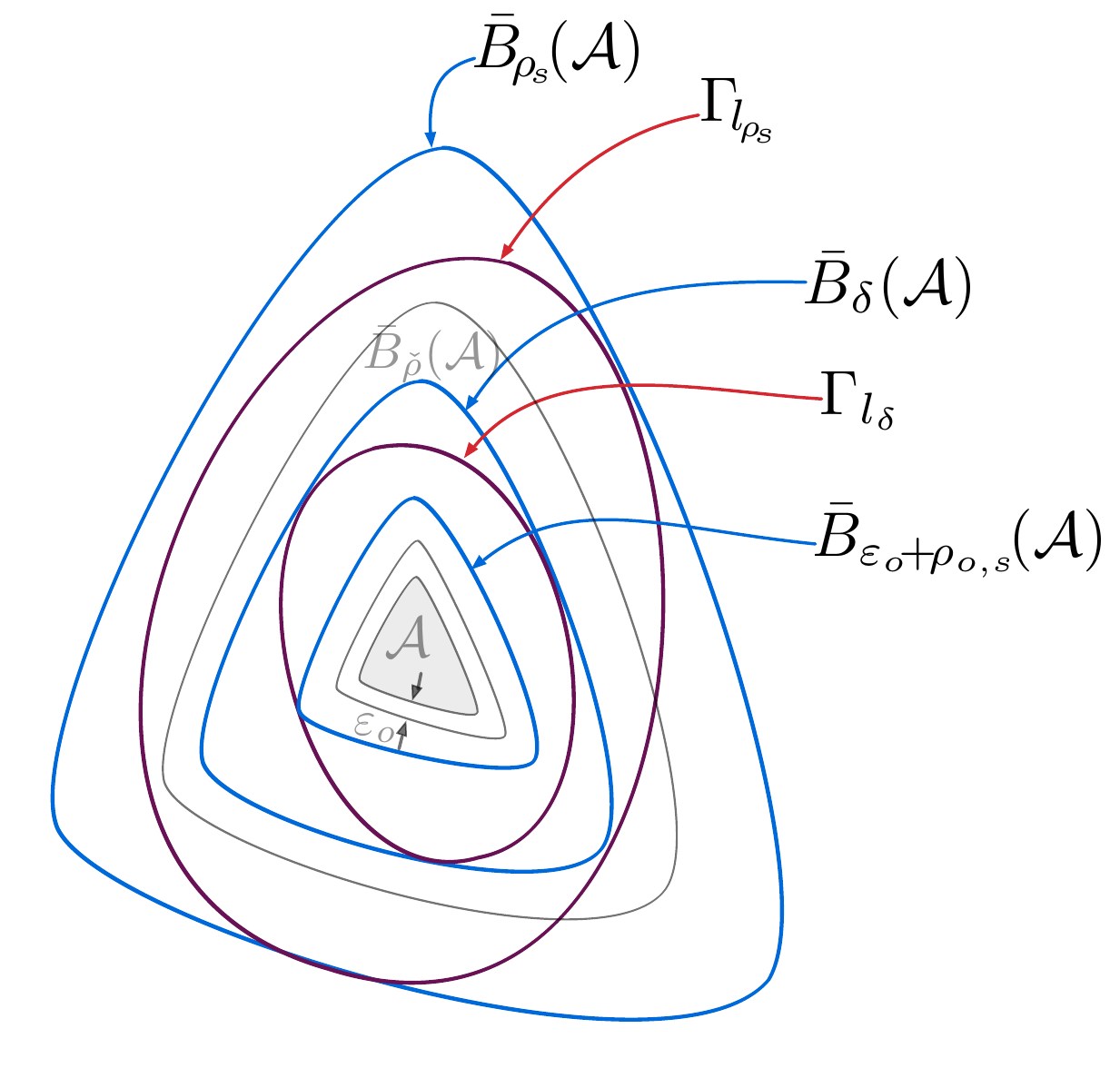}
      \caption{An illustration of Construction \ref{con:backward} used in the proof of practical stability of $\mathcal{A}$ in Theorem \ref{thm:spas}. Whereas Construction \ref{con:inner} starts from $\bar{B}_{\epsilon_o}(\mathcal{A})$ (shown in faded grey), Construction \ref{con:backward} starts from the desired $\bar{B}_{\rho_s}(\mathcal{A})$ to generate the number $\delta$ in Definition \ref{def:stabilityNew2}, and the numbers $\rho_{o,s}$ and $l_{\rho_s}$ that determine the parameter set $P_s$. We have $\rho_{o,s}>0$ because $\rho_s>\check{\rho}$.
      }
      \label{fig:backwardConstructions}
      \vspace*{-.1in}
   \end{figure}

Construction \ref{con:backward} is summarized by the following set containment relationships:
\begin{equation}\label{eq:summaryBackward}
\bar{B}_{\rho_s}(\mathcal{A})\supseteq\Gamma_{l_{\rho_s}}\supseteq\bar{B}_{\delta}(\mathcal{A})\supseteq \Gamma_{l_{\delta}}\supseteq \bar{B}_{\epsilon_o+\rho_{o,s}}(\mathcal{A}),
\end{equation}
which are depicted in Figure \ref{fig:backwardConstructions}.

\begin{rem}
From the preceding constructions, we make the following observations. For $\rho_s=\check{\rho}_s$, Construction \ref{con:backward} coincides with Construction \ref{con:inner}, in the sense that $l_{\rho}=\check{l}$, $\delta=\check{\delta}$ and $l_{\delta}=l_o$. Whenever $\rho_s>\check{\rho}_s$, evidently  $\rho_{o,s}>0$ and $\delta>\check{\delta}$. Moreover, for any $\rho_s\in (\check{\rho}_s,\sigma)$, the construction of $l_{\rho_s}$ in \eqref{eq:lRho} obeys $l_{\rho_s}\in(\check{l},\hat{l})$.
$\diamondsuit$\end{rem}

The practical stability of $\mathcal{A}$ for \eqref{eq:system} follows from Constructions \ref{con:inner} and \ref{con:backward}, and by applying Claim \ref{cl:invarianceNew}. Specifically, generate the number $\check{\rho}_s$ using Construction \ref{con:inner}, and the numbers $\delta$, $\rho_{o,s}$ and $l_{\rho_s}$ using Construction \ref{con:backward}. Then, apply Claim \ref{cl:invarianceNew}, with $\tilde{\rho}=\rho_{o,s}$,  $\tilde{\sigma}=\delta$, and $l=l_{\rho_s}$\footnote{Note that in the proof of Claim \ref{cl:invarianceNew}, Construction \ref{con:outer} is applied with $\tilde{\sigma}=\delta$, generating the number $\hat{l}$, which, according to \eqref{eq:lRho}, coincides with $l_{\rho_s}$.}. This generates the parameter set $P_s$ such that  $\Gamma_{l_{\rho_s}}$ is forward invariant for \eqref{eq:system} whenever $\pi\in P_s$. We observe that Definition \ref{def:stabilityNew2} is satisfied since trajectories initialized inside $\bar{B}_{\delta}(\mathcal{A})\cap\Xi$ never leave $(\Gamma_{l_{\rho_s}}\cap\Xi)\supseteq(\bar{B}_{\delta}(\mathcal{A})\cap\Xi)$, and therefore always remain within $(\bar{B}_{\rho_s}(\mathcal{A})\cap\Xi)\supseteq (\Gamma_{l_{\rho_s}}\cap\Xi)$.

\subsubsection{Semiglobal, Practical Attractivity}\label{sec:attractivity}
Our second task is to show that under the hypotheses of the Theorem, $\mathcal{A}$ is semiglobally, practically attractive for \eqref{eq:system}, in the sense of Definition \ref{def:attractivityNew2}.

We pick $\check{\rho}_a=\epsilon_o$. As in Definition \ref{def:attractivityNew2}, let $\sigma$ and $\rho_a$ be arbitrary, with $\sigma>\rho_a>\check{\rho}_a$.
Use Construction \ref{con:outer}, with $\tilde{\sigma}=\sigma$, to generate the numbers $\hat{\sigma}$ and $\hat{l}$. Generate a parameter set $P_a$ by applying Claim \ref{cl:invarianceNew}, with $\tilde{\rho}=\rho_a-\check{\rho}_a$, $\tilde{\sigma}=\sigma$ and $l=\hat{l}$. Thus, by choosing $\pi\in P_a$, we ensure that $\Gamma_{\hat{l}}$ is forward invariant for \eqref{eq:system}, and that
\begin{align}
& \Delta V(\xi) < -W_{\hat{\sigma},\check{\rho}_a}(\xi), \quad \forall \xi\in \big(\bar{B}_{\hat{\sigma}} \backslash B_{\rho_a}(\mathcal{A})\big)\cap\Xi,\quad \text{and}
\label{eq:DeltaVneg_a}
\\
&W_{\hat{\sigma},\check{\rho}_a}(\xi)>0 \quad \forall \xi\in \big(\bar{B}_{\hat{\sigma}} \backslash B_{\rho_a}(\mathcal{A})\big)\cap\Xi.
\label{eq:Wpos}
\end{align}

In order to show that $\mathcal{A}$ is semiglobally, practically attractive for \eqref{eq:system} on $\bar{B}_{\sigma}(\mathcal{A})$, we will show that $\bar{B}_{\rho_a}(\mathcal{A})$ is uniformly attractive for \eqref{eq:system} on $\bar{B}_{\sigma}(\mathcal{A})$, whenever $\pi\in P_a$. Let $\varepsilon\in\mathbb{R}_{++}$ be arbitrary, but such that $\bar{B}_{\varepsilon}(\bar{B}_{\rho_a}(\mathcal{A}))\cap\Xi\subset \bar{B}_{\sigma}(\mathcal{A})\cap\Xi$. Since $W_{\hat{\sigma},\epsilon_o}(\cdot)$ is continuous, it attains a minimum value on the compact set $\bar{B}_{\hat{\sigma}}(\mathcal{A})\backslash B_{\varepsilon}(\bar{B}_{\rho_a}(\mathcal{A}))$. Let
\begin{align}\label{eq:minW}
\gamma := &\min_{\xi\in\mathbb{R}^n} W_{\hat{\sigma},\check{\rho}_a}(\xi)
\\
&\text{s.t.}\quad \xi\in(\bar{B}_{\hat{\sigma}}(\mathcal{A})\backslash B_{\varepsilon}(\bar{B}_{\rho_a}(\mathcal{A})))\cap\Xi
\nonumber
\end{align}
By \eqref{eq:Wpos}, $W_{\hat{\sigma},\epsilon_o}(\cdot)>0$ on $\bar{B}_{\hat{\sigma}}(\mathcal{A})\backslash B_{\rho_a}(\mathcal{A})$,  and therefore $\gamma$ is strictly positive. It follows from \eqref{eq:DeltaVneg_a} and the fact that $\bar{B}_{\hat{\sigma}}(\mathcal{A})\supset \bar{B}_{\sigma}(\mathcal{A})$, that whenever $\xi(0)\in(\bar{B}_{\sigma}(\mathcal{A})\backslash B_{\varepsilon}(\bar{B}_{\rho_a}(\mathcal{A})))\cap\Xi$, the inequality
\begin{equation}\label{eq:Vev}
V(\xi^+)\leq V(\xi)-\gamma
\end{equation}
holds for all $t\in\mathbb{N}$ for which $\xi(t)\in (\bar{B}_{\hat{\sigma}}(\mathcal{A})\backslash B_{\varepsilon}(\bar{B}_{\rho_a}(\mathcal{A})))\cap\Xi$.
Solving \eqref{eq:Vev}, we obtain that
\begin{equation}
V(\xi(t)) \leq V(\xi(0))-t\gamma,
\end{equation}
which, together with the positive definiteness of $V(\cdot)$, shows that no sequence $(\xi(t))_{t=0}^{\infty}$, generated by \eqref{eq:system} and initialized inside $\bar{B}_{\sigma}(\mathcal{A})\cap\Xi$ can remain in \\ $(\bar{B}_{\hat{\sigma}}(\mathcal{A})\backslash B_{\varepsilon}(\bar{B}_{\rho_a}(\mathcal{A})))\cap\Xi$ forever. The forward invariance of $\Gamma_{\hat{l}}$ implies that no such sequence can leave $\Bgo{\hat{\sigma}}\cap\Xi$, and we therefore conclude that $(\xi(t))_{t=0}^{\infty}$ enters $B_{\varepsilon}(\bar{B}_{\rho_a}(\mathcal{A})))\cap\Xi$ in finitely many iterations, showing that $\bar{B}_{\rho_a}(\mathcal{A})$ is attractive for \eqref{eq:system} on $\bar{B}_{\sigma}(\mathcal{A})$.

Since $\sigma$ and $\rho_a$ are chosen arbitrarily (but such that $\sigma>\rho_a>\check{\rho}_a$), we have shown that $\mathcal{A}$ is semiblobally practically attractive for \eqref{eq:system}.

Having shown that $\mathcal{A}$ is both practically stable and semiglobally, practically attractive for \eqref{eq:system}, we conclude that $\mathcal{A}$ is semiglobally, practically, asymptotically stable for \eqref{eq:system}, and the theorem is proved. $\square$.

\section{Conclusions}\label{sec:conclusions}
We considered a general class of constrained, nonlinear, discrete-time, gain-parametrized systems, and we contributed a tool for the qualitative analysis of such systems. We provided a definition of \emph{semiglobal, practical, asymptotic stability} of compact sets under a given dynamic, and a theorem that characterizes this property in terms of conditions on a Lyapunov-like function. The SPAS theorem we provided does not require the existence of an asymptotically stable attractor for a nominal counterpart to a given system dynamic, and a set having the SPAS property for a given system need not constitute a set of fixed points for that system.

\section*{Acknowledgments}
The author would like to thank Professor Andrew Teel for valuable feedback, and to Professor Naomi Leonard for her support and encouragement.

\newpage
 \bibliographystyle{unsrt}
 \bibliography{ms}
\end{document}